\documentclass[12pt,a4,paper]{article}
\usepackage{amsmath}
\usepackage{amssymb}
\usepackage{bbm}

\def\RR{\mathbb{R}}
\def\CC{\mathbb{C}}
\def\KK{\mathbb{K}}
\def\NN{\mathbb{N}}

\def\sgn{{\rm sgn}\,}

\newtheorem{theorem}{Theorem}

\newtheorem{proposition}[theorem]{Proposition}

\def\be{\begin{enumerate}}
\def\ee{\end{enumerate}}

\begin{document}

\begin{center} \Large \textbf{A remark on the rigidity of a property characterizing the Fourier transform} \end{center}

\begin{center} Hermann K\"onig, Vitali Milman \end{center}
\vspace{1,5ex}

Consider the classical Fourier transform on $\RR^n$ given by
$$\mathcal{F} f(x) = \int_{\RR^n} \exp(-2 \pi i \langle x , y \rangle ) \ f(y) \ dy $$
e.g. on the Schwartz space $\mathcal{S}(\RR^n)$ of rapidly decreasing smooth functions $f: \RR^n \to \CC$. As well-known, $\mathcal{F}$ acts bijectively on $\mathcal{S}(\RR^n)$ and exchanges products with convolutions. In a series of papers by Alesker, Artstein-Avidan, Faifman and Milman \cite{AAM}, \cite{AAFM} and \cite{AFM} it was shown that these properties essentially characterize the Fourier transform. The final result in \cite{AFM} is: Any bijective transformation $T : \mathcal{S}(\RR^n) \to \mathcal{S}(\RR^n)$ satisfying $T(f \star g) = T f \cdot Tg$ for all $f, g \in \mathcal{S}(\RR^n)$ is just a slight modification of the Fourier transform: there exists a diffeomorphism $\omega: \RR^n \to \RR^n$ such that either $Tf = (\mathcal{F} f) \circ \omega$ for all $f \in \mathcal{S}(\RR^n)$ or $Tf = \overline{(\mathcal{F} f) \circ \omega}$ for all $f \in \mathcal{S}(\RR^n)$. If in addition $T(f \cdot g) = Tf \star Tg$ holds for all $f, g \in \mathcal{S}(\RR^n)$, the diffeomorphism is given by a linear map $\omega = A$, see \cite {AFM}. Note that $T$ is not assumed to be linear or continuous. Nevertheless the real linearity and the continuity are a consequence of the result. \\

The idea in the proof of \cite{AFM} is to consider the map $S = T \circ \mathcal{F}^{-1}: \mathcal{S}(\RR^n) \to \mathcal{S}(\RR^n)$, which is multiplicative, $S(F \cdot G) = SF \cdot SG$ for all $F, G \in \mathcal{S}(\RR^n)$. The paper \cite{AFM} characterizes bijective multiplicative maps between suitable function spaces like $\mathcal{S}(\RR^n)$ as having the form $SF = F \circ \omega$ or $SF = \overline{F \circ \omega}$ for a diffeomorphism $\omega$. Already in 1949 Milgram [Mi] had studied  bijective multiplicative operators $S : C(M) \to C(M)$ on real-valued spaces of continuous functions $C(M)$ on a finite dimensional manifold $M$. They have the form $Sf(x) = |f(\omega(x))|^{p(x)} \ \sgn f(\omega(x))$, where $p:M \to \RR_{>0}$ is a suitable continuous function and $\omega : M \to M$ is a homeomorphism. This result was extended to spaces of smooth functions $C^k(M)$ on a differentiable $C^k$-manifold $M$ by Mr\v{c}un and \v{S}emrl \cite{MS} for $k \in \NN$. It was extended to $C^\infty(M)$-functions and other function algebras like the Schwartz space $\mathcal{S}(\RR^n)$ and to spaces of complex-valued functions in \cite{AFM}, see also the earlier paper \cite{AAFM}. The proof in \cite{AFM} also works in the case that the constant function $1$ does not belong to the space, as in the case of $\mathcal{S}(\RR^n)$. As indicated, the result for $\mathcal{S}(\RR^n)$ is useful to characterize the Fourier transform as essentially a bijective map $T$ on $S(\RR^n)$ which maps convolutions to products and vice-versa, yielding a characterization up to diffeomorphisms $\omega$ and complex conjugation, $Tf = (\mathcal{F} f) \circ \omega$ or $Tf= \overline{(\mathcal{F} f) \circ \omega}$, see \cite{AFM}. \\

We examine the question whether these characterizations are rigid under perturbations. Again we do not assume linearity or continuity of the operators $T$ under consideration. As for the spaces on which $T$ operates, let us assume that $M$ is a finite dimensional $C^k$-manifold, $k \in \NN \cup \{0\} \cup \{\infty \}$ (for $k=0$ just a manifold) and that
$$E \subset C^k(M) := \{f : M \to \KK \ | \ f \text{ \ k-times continuously differentiable} \}$$
is a linear subspace, where $\KK \in \{\RR, \CC \}$. We let $C(M) = C^0(M)$. We call $E$ {\it rich} if $f, g \in E$ always implies $f \cdot g \in E$ and if for all $x \in M$ and all $r>0$ there is $h \in E$ with $|h(x)| > r$. The following rigidity result holds for bijective multiplicative operators. \\

\begin{proposition}
Let $k \in \NN \cup \{0\} \cup \{\infty \}$, $M$ be a finite dimensional $C^k$-manifold and $E \subset C^k(M)$ be a rich subspace. Let $a : M \to \RR_{>0}$ be a function and assume that $T : E \to E$ is a bijective operator such that for all $f, g \in E$ and all $x \in M$
\begin{equation}\label{eq1}
| T(f \cdot g)(x) - Tf(x) \cdot Tg(x) | \le a(x)
\end{equation}
holds. Then
\begin{equation}\label{eq2}
T(f \cdot g) = Tf \cdot Tg \text{ \ for all \ } f, g \in E \ .
\end{equation}
\end{proposition}

Thus almost multiplicativity implies multiplicativity: The operation is rigid. Note that $T$ is not assumed to be linear; otherwise \eqref{eq1} would directly imply \eqref{eq2}. Since we fix a point $x$ with \eqref{eq1} and prove \eqref{eq2} for this point $x$, \eqref{eq1} implies \eqref{eq2} for any subset $A$ of $M$. Choosing $A$ to be a dense subset of $M$ yields that, if condition \eqref{eq1} holds for all $x \in A$, \eqref{eq2} will follow for all $f, g \in E$ on $\bar{A} = M$, since $T(f \cdot g)$ and $Tf \cdot Tg$ are continuous. \\

{\bf Proof.}
Let $g \in E$ and $x \in M$ be fixed. Since $E$ is rich, for any $r>0$ there is $k_r \in E$ with $|k_r(x)| > r$. Since $T$ is bijective, there is $h_r \in E$ with $|(T h_r)(x)| = |k_r(x)| > r$. The assumption implies
$$| \frac{T(g \cdot h_r)(x)}{(T h_r)(x)} - Tg(x) | \le \frac {a(x)}{|(T h_r)(x)|} \le \frac{a(x)} r \ . $$
Hence $Tg(x) = \lim_{r \to \infty} \frac{T(g \cdot h_r)(x)}{(T h_r)(x)}$. \\
Next, consider $x \in M$ and $f, g \in E$ as well as $h_r \in E$ with $|(T h_r)(x)| > r$. Then
$$| T(f \cdot g \cdot h_r)(x) - T(f \cdot g)(x) \cdot (T h_r)(x) | \le a(x) \ , $$
$$| T(f \cdot g \cdot h_r)(x) - Tf(x) \cdot T(g \cdot h_r)(x) | \le a(x) \ . $$
Therefore
$$| Tf(x) \cdot T(g \cdot h_r)(x) - T(f \cdot g)(x) \cdot (T h_r)(x) | \le 2 a(x) \ , $$
$$| Tf(x) \cdot \frac{T(g \cdot h_r)(x)}{(T h_r)(x)} - T(f \cdot g)(x) | \le \frac{2 a(x)}{|(T h_r)(x)|} \le \frac{2 a(x)} r \ . $$
Hence $T(f \cdot g)(x) = \lim_{r \to \infty}  Tf(x) \cdot \frac{T(g \cdot h_r)(x)}{(T h_r)(x)} = Tf(x) \cdot Tg(x)$ and $T$ is multiplicative.
\hfill $\Box$   \\

For $k \in \NN \cup {0} \cup {\infty}$, $E = C^k(M)$ and $E = \mathcal{S}(\RR^n)$ are rich subspaces. For $k \in \NN \cup \{\infty\}$, by \cite{MS} and \cite{AFM}, equation \eqref{eq2} implies for bijective operators $T: E \to E$ that $Tf = f \circ \omega$ for all $f \in E$ or $Tf = \overline{f \circ \omega}$ for all $f \in E$ (in the complex case) where $\omega : M \to M$ is a suitable $C^k$-diffeomorphism, in the case of $E = \mathcal{S}(\RR^n)$ a $C^\infty$-diffeomorphism. \\ 
For continuous functions on $M$, $k=0$, by \cite{M}, we have that $Tf(x) = |f \circ \omega(x)|^{p(x)} \ \sgn(f \circ \omega(x))$ for all $f \in E$ (in the real case), where $\omega: M \to M$ is a homeomorphism and $p: M \to \RR_{>0}$ is a continuous function. \\

Let us call a subspace $E \subset C^k(\RR^n)$, $k \in \NN \cup {0} \cup {\infty}$ {\it convolution-stable}, if for all $f, g \in E$ we have $f \star g \in E$. The following rigidity result holds for the convolution operator. \\

\begin{proposition}\label{prop2}
Let $n \in \NN$, $k \in \NN \cup \{0\} \cup \{\infty \}$, $E \subset C^k(\RR^n)$ be a convolution-stable subspace and $F \subset C^k(\RR^n)$ be a rich subspace. Let $a : \RR^n \to \RR_{>0}$ be a function and assume that $T : E \to F$ is a bijective operator such that for all $f, g \in E$ and all $x \in \RR^n$
\begin{equation}\label{eq3}
| T(f \star g)(x) - Tf(x) \cdot Tg(x) | \le a(x)
\end{equation}
holds. Then
\begin{equation}\label{eq4}
T(f \star g) = Tf \cdot Tg \text{ \ for all \ } f, g \in E \ .
\end{equation}
\end{proposition}

Again, $T$ is assumed neither real-linear nor continuous. \\

{\bf Proof.}
Let $g \in E$, $x \in \RR^n$ and $r>0$. Since $F$ is rich, there is $k_r \in F$ with $|k_r(x)| > r$. Since $T: E \to F$ is bijective, there is $h_r \in E$ with
$|(T h_r)(x)| = |k_r(x)| > r$. Hence
$$| \frac{T(g \star h_r)(x)}{(T h_r)(x)} - Tg(x) | \le \frac {a(x)}{|(T h_r)(x)|} \le \frac{a(x)} r \ . $$
Hence $Tg(x) = \lim_{r \to \infty} \frac{T(g \star h_r)(x)}{(T h_r)(x)}$. \\
Similarly as above, this implies for $f, g \in E$ by considering $T(f \star g \star h_r)$ that
$$| Tf(x) \cdot T(g \star h_r)(x) - T(f \star g)(x) \cdot (T h_r)(x) | \le 2 a(x)  $$
and
$$T(f \star g)(x) = \lim_{r \to \infty}  Tf(x) \cdot \frac{T(g \star h_r)(x)}{(T h_r)(x)} = Tf(x) \cdot Tg(x) \ , $$
thus $T(f \star g) = Tf \cdot Tg$.
\hfill $\Box$   \\

The proof shows that if \eqref{eq3} holds for one point $x$, \eqref{eq4} will also hold for $x$. Hence if \eqref{eq3} holds for all $x$ in a dense subset $A$ of $\RR^n$, \eqref{eq4} will be true on $\bar{A} = \RR^n$. Note that $T(f \star g)$ and $Tf \cdot Tg$ are continuous by assumption. \\

We may apply this to $E = F = \mathcal{S}(\RR^n)$. Then equation \eqref{eq4} implies that there exists a $C^\infty$-diffeomorphism $\omega : \RR^n \to \RR^n$ that either $Tf = (\mathcal{F} f) \circ \omega$ for all $f \in E$ or $Tf = \overline{(\mathcal{F} f) \circ \omega}$ for all $f \in E$. \\
In the continuous case of Proposition \ref{prop2}, $k=0$, there exists additionally a continuous function $p : \RR^n \to \RR_{>0}$ such that $Tf(x) = |(\mathcal{F} f) \circ \omega (x)|^{p(x)} \ \sgn((\mathcal{F} f) \circ \omega(x))$ for all $f \in E$ and $x \in \RR^n$. This follows by applying \eqref{eq2} to the multiplicative map $S := T \mathcal{F}^{-1} : F \to F$, see Theorems 1.3 and 1.10 of \cite{AFM}. \\
Another application of \eqref{eq3} and \eqref{eq4} is given by the choice $E = C^\infty_c(\RR^n,\CC)$, the complex-valued $C^\infty$-functions with compact support and $F = \mathcal{F}(E) = PW(\RR^n,\CC)$, the Paley-Wiener space of functions $f : \RR^n \to \CC$ decreasing polynomially on $\RR^n$ and being extendable to $F : \CC^n \to \CC$ with $|F(z)| \le A \exp(B |z|)$ for all $z \in \CC^n$, where $A, B >0$ are constants. Then by Theorem 1.5 of \cite{AFM} $T(f \star g) = Tg \cdot Tg$ also implies that either $Tf = (\mathcal{F} f) \circ \omega$ for all $f \in E$ or $Tf = \overline{(\mathcal{F} f) \circ \omega}$ for all $f \in E$.

\begin{proposition}
Let $n \in \NN$, $E = \mathcal{S}(\RR^n)$ and $a > 0$ be a constant. Assume that $T : E \to E$ is a bijective operator satisfying
\begin{equation}\label{eq5}
| T(f \cdot g)(x) - (Tf \star Tg)(x) | \le a \
\end{equation}
for all $f, g \in E$ and all $x \in \RR^n$.
Then
\begin{equation}\label{eq6}
T(f \cdot g) = Tf \star Tg
\end{equation}
for all $f, g \in E$.
\end{proposition}

{\bf Proof.}
Let $f \in E$ and $x \in \RR^n$ be fixed. Choose an approximation of the $\delta$-distribution in $F$, i.e. a sequence of $L_1(\RR^n)$-functions
$(l_r)_{r \in\NN}$ such that for all $h \in F$ we have \newline $\lim_{r \to \infty} (h \star l_r)(x) = h(x)$. Since $T$ is bijective, we may choose $\psi_r \in E$ such that $T(\psi_r) = r \cdot l_r$. By \eqref{eq5}
$$ | T(f \cdot \psi_r)(x) - (Tf \star T \psi_r)(x) | \le a \ , $$
$$ \left| \frac{T(f \cdot \psi_r)(x)} r - (Tf \star l_r)(x) \right| \le \frac a r \ , $$
so that
\begin{equation}\label{eq7}
\lim_{r \to \infty} \frac{T(f \cdot \psi_r)(x)} r = \lim_{r \to \infty} (Tf \star l_r)(x) = Tf(x) \ .
\end{equation}
For $f, g \in E$ we have by \eqref{eq5} for all $y \in \RR^n$
$$ | T(f \cdot g \cdot \psi_r)(y) - (Tf \star T(g \cdot \psi_r))(y) | \le a $$
and
$$T(g \cdot \psi_r)(y) - (Tg \star T \psi_r)(y) =: \phi(y) \ , |\phi(y)| \le a \ , $$
with $\phi \in \mathcal{S}(\RR^n)$ depending on $g$ and $\psi_r$. Hence
\begin{equation}\label{eq8}
\left| \frac{T(f \cdot g \cdot \psi_r)(x)} r - \frac{(Tf \star Tg \star  T \psi_r)(x)} r - \frac{(Tf \star \phi)(x)} r \right| \le \frac a r \ .
\end{equation}
Note that $(Tf \star \phi)(x)$ is well defined and bounded independent of $g$ and $\psi_r$, since $Tf \in \mathcal{S}(\RR^n)$ and $\phi \in \mathcal{S}(\RR^n)$ is uniformly bounded by $a$. By \eqref{eq7}
$$\lim_{r \to \infty} \frac{T(f \cdot g \cdot \psi_r)(x)} r = T(f \cdot g)(x) \ , $$
and
$$\lim_{r \to \infty} \frac{(Tf \star Tg \star T \psi_r)(x)} r = \lim_{r \to \infty} (Tf \star Tg \star l_r)(x) = (Tf \star Tg)(x) \ . $$
Therefore \eqref{eq8} implies $T(f \cdot g)(x) = (Tf \star Tg)(x)$ for all $x \in \RR^n$.
\hfill $\Box$   \\

Consider the operator $S = \mathcal{F}^{-1} T : E \to E$. By \eqref{eq6} $S$ is multiplicative, and Theorem 1.3 of \cite{AFM} implies that there is a $C^\infty$-diffeomorphism $\omega: \RR^n \to \RR^n$ such that either $Tf = \mathcal{F}(f \circ \omega)$ for all $f \in E$ or $Tf = \overline{\mathcal{F}(f \circ \omega)}$ for all $f \in E$. \\ 
For $k=0$ we get for real-valued functions $Tf(x) = |\mathcal{F}(f \circ \omega)(x)|^{p(x)} \ \sgn(\mathcal{F}(f \circ \omega)(x))$ for all $f \in E$ and $x \in \RR^n$, where $p : \RR^n \to \RR_{>0}$ is a suitable continuous function. \\

We also prove a rigidity result for the chain rule.

\begin{proposition}
Let $a : \RR \to \RR_{>0}$ be a function. Assume that $T : C^1(\RR) \to C(\RR)$ satisfies
\begin{equation}\label{eq9}
| T(f \circ g)(x) - (Tf)(g(x)) \cdot Tg(x) | \le a(x)
\end{equation}
for all $f, g \in C^1(\RR)$ and all $x \in \RR$. Suppose that for all $x \in \RR$ and $r>0$ there is a bounded function $h_r \in C^1(\RR)$ such that $h_r(x) = x$ and $|(Th_r)(x)| > r$. Then for all $f, g \in C^1(\RR)$
\begin{equation}\label{eq10}
T(f \circ g) = (Tf)\circ g \cdot Tg \ .
\end{equation}
\end{proposition}

Again, $T$ is not assumed to be continuous or linear. If \eqref{eq9} holds for all $x$ in a dense subset of $\RR$, \eqref{eq10} will follow on $\RR$. \\

{\bf Proof.}
Let $x \in \RR$ and $r>0$. By assumption there is $h_r \in C^1(\RR)$ with $h_r(x)=x$ and $|(T h_r)(x)| > r$. By \eqref{eq9} for all $f \in C^1(\RR)$
$$| T(f \circ h_r)(x) -(Tf)(h_r(x)) \cdot (T h_r)(x) | = | T(f \circ h_r)(x) -(Tf)(x) \cdot (T h_r)(x) | \le a(x) \ , $$
$$| \frac{T(f \circ h_r)(x)}{(T h_r)(x)} - (Tf)(x) | \le \frac{a(x)}{|(T h_r)(x)|} \le \frac{a(x)} r \ , $$
so that $(Tf)(x) =\lim_{r \to \infty} \frac{T(f \circ h_r)(x)}{(T h_r)(x)}$. Considering for $f, g \in C^1(\RR)$ the map $f \circ g \circ h_r$, we have by \eqref{eq9}
\begin{align*}
T(f \circ g \circ h_r)(x) & = (Tf)((g \circ h_r)(x)) \cdot T(g \circ h_r)(x) + \psi_1(x) \\
& = (Tf)(g(x)) \cdot T(g \circ h_r)(x) + \psi_1(x) \ ,
\end{align*}
where $\psi_1$ depends on $f$, $g$ and $h_r$, but satisfies $|\psi_1(x)| \le a(x)$. Further
\begin{align*}
T(g \circ h_r)(x) & = (Tg)(h_r(x)) \cdot (T h_r)(x) + \psi_2(x) \\
& = (Tg)(x) \cdot (T h_r)(x) + \psi_2(x) \ ,
\end{align*}
again $\psi_2$ depending on $g$ and $h_r$, but with $|\psi_2(x)| \le a(x)$. We find
\begin{align*}
T(f \circ g)(x) & = \lim_{r \to \infty} \frac{T(f \circ g \circ h_r)(x)}{(T h_r)(x)} = \lim_{r \to \infty} \left( (Tf)(g(x)) \cdot \frac{T(g \circ h_r)(x)}{(T h_r)(x)} + \frac{\psi_1(x)}{(T h_r)(x)} \right) \\
& = (Tf)(g(x)) \cdot (Tg)(x) + \lim_{r \to \infty} \left( \frac{\psi_1(x)}{(T h_r)(x)} + \frac{(Tf)(g(x)) \psi_2(x)}{(T h_r)(x)} \right) \\
& = (Tf)(g(x)) \cdot (Tg)(x) \ .
\end{align*}
Therefore $T(f \circ g) = (Tf) \circ g \cdot Tg$ follows.
\hfill $\Box$   \\

By \cite{AKM}, see also \cite{KM}, equation \eqref{eq10} implies that there is $p \ge 0$ and there is a continuous positive function $H \in C(\RR)$ such that we have either $Tf = \frac{H \circ f} H |f'|^p$ for all $f \in C^1(\RR)$ or in the case $p>0$ that $Tf = \frac{H \circ f} H |f'|^p \ \sgn(f')$ for all $f \in C^1(\RR)$. \\

In the case of the Leibniz rule, we have the following rigidity result. \\

\begin{proposition}
Suppose $T : C^1(\RR) \to C(\RR)$ is an operator and $a : \RR \to \RR_{>0}$ is a function such that for all $f \in C^1(\RR)$ and $x \in \RR$ we have
\begin{equation}\label{eq11}
| T(f \cdot g)(x) - Tf(x) \cdot g(x) - f(x) \cdot Tg(x) | \le a(x) \ .
\end{equation}
Assume for all $x \in \RR$ there are functions $h_n \in C^1(\RR)$ such that $\lim_{n \to \infty} |h_n(x)| = \infty$ and $\lim_{n \to \infty} \frac{|Th_n(x)|}{|h_n(x)|} = 0$. Then for all $f, g \in C^1(\RR)$
\begin{equation}\label{eq12}
T(f \cdot g) = Tf \cdot g + f \cdot Tg \ .
\end{equation}
\end{proposition}

{\bf Proof.}
Let $x \in \RR$ and $h_n$ be as in the statement of the proposition. Then \eqref{eq11} implies
$$\left| \frac{T(f \cdot h_n)(x)}{h_n(x)} - Tf(x)  - f(x) \cdot \frac{Th_n(x)}{h_n(x)} \right| \le \frac{a(x)}{|h_n(x)|} \ , $$
so that $ \lim_{n \to \infty} \frac{T(f \cdot h_n)(x)}{h_n(x)} = Tf(x)$. This implies for $f, g \in C^1(\RR)$ that
$$\left| \frac{T(f \cdot g \cdot h_n)(x)}{h_n(x)} - f(x) \cdot \frac{T(g \cdot h_n)(x)}{h_n(x)} - g(x) \cdot Tf(x) \right| \le \frac{a(x)}{|h_n(x)|} \ , $$
which implies that $T(f \cdot g) = Tf \cdot g + f \cdot Tg$.
\hfill $\Box$   \\

Equation \eqref{eq12} implies by \cite{KM} or \cite{KM1} that there are continuous functions $c, d : \RR \to \RR$ such that $T$ has the form
$$Tf(x) = c(x) \ f'(x) + d(x) \ f(x) \ \ln |f(x)| $$
for all $f \in C^1(\RR)$,

\vspace{1cm}

H. K\"onig\\
Mathematisches Seminar\\
Universit\"at Kiel\\
24098 Kiel, GERMANY  \\

V. Milman\\
School of Mathematical Sciences\\
Tel Aviv University\\
Ramat Aviv, Tel Aviv 69978, ISRAEL


\begin{thebibliography}{99}
\bibitem[KM1]{KM1} H. K\"onig, V. Milman; Characterizing the derivative and the entropy function by the Leibniz rule, with an appendix by D. Faifman,  J. Funct. Anal. 261 (2011), 1325-1344.
\end{thebibliography}
\end{document}